\documentclass{article}
\usepackage{graphicx} 

\usepackage{fancyhdr}
\usepackage{isomath}
\usepackage{mathtools} 
\usepackage{amsbsy}
\usepackage{amssymb}
\usepackage{amscd,amsfonts}
\usepackage{bigints}
\usepackage{graphicx}
\usepackage{verbatim}
\usepackage{euscript}
\usepackage{alltt}
\usepackage{stmaryrd}
\usepackage{relsize}
\usepackage{enumerate}
\usepackage{url}
\usepackage[stable]{footmisc}
\usepackage{breakurl}
\usepackage{hyperref}
\usepackage{comment}
\usepackage[numbers]{natbib}
\usepackage[font=small,labelfont=bf]{caption}
\usepackage[font=small,labelfont=bf]{subcaption}
\usepackage{float}
\usepackage{mwe}
\usepackage{algorithm}
\usepackage{algpseudocode}
\usepackage{xcolor}
\usepackage[super]{nth}
\usepackage{soul}
\usepackage{centernot}
\DeclareGraphicsExtensions{.eps,.pdf}
\usepackage{amsmath}
\usepackage{amsbsy}
\usepackage{amssymb}
\usepackage{amscd}
\usepackage{amsfonts}

\newcommand{\R}{\mathbb R}

\newcommand{\beq}{\begin{equation}}
\newcommand{\eeq}{\end{equation}}
\newcommand{\beqs}{\begin{eqnarray}}
\newcommand{\eeqs}{\end{eqnarray}}
\newcommand{\beql}{\begin{equation} \label}
\newcommand{\half}{\frac{1}{2}}

\newcommand{\calA}{{\cal A}}
\newcommand{\calB}{{\cal B}}
\newcommand{\calC}{{\cal C}}

\newcommand{\calF}{{\cal F}}

\newcommand{\calN}{{\cal N}}

\newcommand{\p}{\partial}

\usepackage[margin=1in]{geometry}

\newcommand{\dee}{\mathcal{D}}
\newcommand{\scl}{\mathcal{L}}

\date{}
\begin{document}
\title{A Hidden Convexity in Continuum Mechanics, with application to classical, continuous-time, rate-(in)dependent plasticity}

\author{Amit Acharya\thanks{Department of Civil \& Environmental Engineering, and Center for Nonlinear Analysis, Carnegie Mellon University, Pittsburgh, PA 15213, email: acharyaamit@cmu.edu.}}

\maketitle
\begin{abstract}
\noindent A methodology for defining variational principles for a class of PDE models from continuum mechanics is demonstrated, and some of its features explored. The scheme is applied to quasi-static and dynamic models of rate-independent and rate-dependent, single crystal plasticity at finite deformation.

\end{abstract}

\section{Introduction}
In this paper we explore a strategy for designing variational principles for a significant class of static and dynamical models from continuum mechanics, naturally stated as systems of partial differential equations (PDE). The models can be dissipative or conservative. Action functionals are designed, whose Euler-Lagrange equations recover the primal PDE system and side conditions in a well-defined sense. The essential ideas behind the approach may be understood from \cite[Sec.~2]{action_3}, \cite[Sec.~7]{action_2}, and \cite[Sec.~6.1]{acharya2022action}. The variational principles govern a dual set of fields corresponding to the primal ones of the continuum mechanical model, and the scheme provides a mapping to recover the primal fields with the guarantee that the latter are weak solutions of the primal model. The Lagrangian of the dual variational problem is convex (with a trivial sign change), and therefore existence of a minimizer appears to rest on only the coercivity of the dual functional. Correspondingly, the Euler-Lagrange equations of the dual functional are shown to possess a local degenerate ellipticity, regardless of the properties of the primal system. In the context of solving the primal PDE system, these features are crucially enabled by the `free' choice of a (family of) potential(s) in the primal variables that may be interpreted as defining a `target' whose integral is to be extremized, subject to the primal PDE system as constraints. The dual fields then are simply the Lagrange multipliers of the formulation, and since the target is free to choose, one chooses it to have as strongly positive-definite an Hessian as needed to dominate the non-monotonicity of the constraint equations. Everything said above is of recent origin \cite{acharya2022action, action_2,action_3} and mathematically formal, but potentially useful, as borne out by encouraging results in computational implementations of model problems \cite{KA1}-\cite[Sec.~6]{a_arora_thesis}-\cite{sga,KA2} involving a range of linear and nonlinear, ODE and PDE, time-dependent and independent problems related to continuum mechanics (linear transport, heat equation, Euler's equations for a rigid body, double-well elastostatics in 1-d, inviscid Burgers in conservation and Hamilton-Jacobi form, the inverse problem of a liquid crystal membrane attaining a prescribed shape, constrained to meet a prescribed principal stretch field), and all approximated by the simplest Galerkin discretization (in these first instances) for solving boundary value problems in domains in space(-time). The work \cite{sga} also contains new, rigorous results on the existence of a variational dual solution for the Saint-Venant Kirchhoff model of nonlinear elasticity whose energy density is not quasiconvex (and hence existence of variational solutions by any other means is not known).

The essential idea behind the approach may be simply understood as follows: Suppose one is interested in solving the generally nonlinear set of equations $G(U) = 0, U \in \R^n, G:\R^n \to \R^m$. The goal is to convert this question of solving a system of equations to an optimization problem, preferably a convex one. To do this, consider, for the moment, an objective function $U \mapsto H(U) \in \R$ and the governing equations for the critical points of the optimization of $H$ subject to the constraints $G(U) = 0$. These equations are, for Lagrange multipliers $\lambda \in \R^m$,
\begin{subequations}
\begin{align}
H(U) - \lambda \cdot G(U) & =:  \scl_H(U,\lambda)\notag\\
   \nabla H(U) - \lambda \cdot \nabla G(U) & =  \p_U \scl_H (U, \lambda)  = 0 \label{eq:intro_LM}\\
    G(U) & = \p_\lambda \scl_H (U, \lambda) = 0 \notag,
\end{align}
\end{subequations}
and the interest would be in solving for a pair $(U,\lambda)$ that satisfies the above set of equations. One now makes the observation that in this formulation of the question of solving $G(U) = 0$ for a $U$, given the function $G$, the choice of $H$ is quite arbitrary and hence, one approach to solving the above set of equations is to solve \eqref{eq:intro_LM} for $U$ in terms of $\lambda$ by making an appropriate choice of $H$ to construct a function $U^{(H)}:\mathcal{O} \to \R^n, \mathcal{O} \subset \R^m$ which satisfies \eqref{eq:intro_LM} for all $\lambda \in \mathcal{O}$, and then look for a solution of $\lambda \in \mathcal{O}$ that satisfies $G\left(U^{(H)}(\lambda)\right) = 0$. This approach defines an adapted change of variables for the problem and ensures that the number of the (dual) variables to be solved for is maintained at the number of the primal equations (the idea carries over to fields as well); more details along these lines are provided in Appendix \ref{app:rev_duality}. Clearly, to the extent that the convexity of the function $H$ dominates any noncovexity of $\lambda \cdot G$ in $U$ for each fixed $\lambda \in \mathcal{O}$, the easier it is to construct the function $U^{(H)}$. Further discussion on the choice of the function $H$ and its effect on the scheme can be found in the preamble of \cite[Sec.~5]{sga}.

An important property of the scheme is that the function $\widehat{S}_H: \R^n \times \R^m \to \R$ defined by $\widehat{S}_H(U,\lambda) = H(U) - \lambda \cdot G(U) = \scl_H(U, \lambda)$ is affine in $\lambda$ and hence concave as well, and therefore $S_H(\lambda) := \inf_U \widehat{S}_H(U, \lambda)$ is concave (but not necessarily strictly so) so that, roughly speaking, obtaining a critical point by the procedure above may be thought of as a problem of concave maximization (or convex minimization), regardless of the nature of the nonlinearities of $G$. 

When $G = 0$ is a system of ordinary or partial differential equations (could involve inequalities too, as in this paper), $\widehat{S}_H, S$ become functionals with $\scl_H$ the Lagrangian for $\widehat{S}_H$ and the critical point equations involve variational derivatives. For boundary-value-problems, all `primal' boundary conditions (whether Dirichlet or Neumann) become natural boundary conditions on the dual side, in quite an `obvious' manner. For primal initial-(boundary)-value-problems an added insight is required; the dual Euler-Lagrange equations are second-order in time and require final-time boundary conditions and it is known that initial value problems do not admit such boundary conditions. It turns out that the dual problem admits (rather arbitrary) Dirichlet boundary conditions on the \emph{dual} fields without affecting the primal solution obtained through the mapping $U^{(H)}$, say when the primal problem has uniqueness of solutions. This has been demonstrated in \cite[Sec.~7]{action_2}-\cite{KA1, KA2}-\cite[Sec.~6.3.4]{sga}.

As examples of this overall approach, in this paper we apply the strategy to develop action principles for classical rate-dependent and independent, dynamic and quasi-static, single crystal plasticity theories without restriction to rate problems, time discretization, energy minimizing paths, associated plasticity, hardening matrix derived from an energy potential, treating plastic slip as an energetic state variable, the existence of a dissipation potential or even a free energy function. Variational principles for plasticity is a subject with a substantial body of work, e.g.~\cite{hill1958general, hill1959some, hill1979aspects, suquet1988discontinuities, strang1979family}-\cite[and earlier references therein]{petryk2003incremental, petryk2020quasi}-\cite{ortiz1989symmetry, ortiz1999nonconvex, ortiz1999variational, carstensen2002non, mielke2015rate, maso2006quasistatic, dal2011quasistatic}, with a detailed review presented in \cite{petryk2020quasi}. Our work is complementary to these points of view, and presents a different formalism that exploits the `free' choice of an added target function in the primal variables with strong convexity properties, especially in providing a somewhat unified point of view in dealing with quasi-static and dynamic problems without utilizing time-discretization. The possibility of making such a  choice in aiding the solution to problems has the flavor of the use of the `linear comparison solid' related to the literature on effective properties, cf., \cite{castaneda1997nonlinear, castaneda1999variational}.

Rigorous results on a weak formulation and existence of solutions of the governing PDE of plasticity were first provided in the seminal work \cite[and earlier references]{suquet1988discontinuities} and subsequent works , e.g.~\cite{han2012plasticity}. As is well-understood by experts but simply to avoid possible confusion, we explicitly note that the existence of a variational principle for a set of (O)PDE is a different question than that of posing a variational (weak) statement for that set of equations. 

An outline of the paper is as follows: in Sec.~\ref{sec:dual_cont_mech} we present the dual formulation. In Sec.~\ref{sec:deg_ellipt} a computation is presented to motivate the degenerate ellipticity of the dual problem. Secs.~\ref{sec:rd_sc} and \ref{sec:ri_sc} contain the algorithmic steps to implement the scheme on the theories of rate-dependent and rate-independent single crystal plasticity, respectively. Sec.~\ref{sec:concl} contains some concluding remarks.

A few words on notation: \emph{except for Sec.~\ref{sec:ri_sc}},  we always use the summation convention on ranges of indices, and the placement of indices as super or subscripts has no special siginifcance. The use of direct notation would have required too many definitions to be put in place for many of the explicit calculations - hence, despite their clumsy appearance, I have chosen to explicitly write out the computations - it is hoped that this avoids any ambiguity. Also, whenever a function is declared as capable of being defined arbitrarily, such arbitrariness is assumed to be restricted by natural smoothness requirements for the problem context to make sense.

We mention at the very outset that from the point of view of this paper, the variational principles developed are purely mathematical devices with their sole justification resting on contributing to solution strategies for the primal system of physical equations involved. Thus, the whole burden of physical modeling rests on the development of the primal system, which is considered a `given' in this work. This paper is not concerned with the quality of physical modeling of plasticity with the considered models, or their connection to the micromechanics of the phenomenon. Those are separate concerns dealt with in \cite{arora2020dislocation, arora2020unification, arora2020finite, arora2022mechanics, arora2023interface} - in fact, this more sophisticated model (which incorporates many of the features of the classical theory), which, among other things, is a first example of a setting in continuum solid mechanics where non-singular, finite deformation elastic fields of arbitrary dislocation distributions can be calculated, served as the primary motivation for the development of the formalism presented here, as discussed in \cite{action_3}.
\section{A dual formulation for models from continuum mechanics}\label{sec:dual_cont_mech}
This paper started out with the specific goal of demonstrating some variational principles for the equations of plasticity theory. However, I soon realized that the main ideas were most efficiently conveyed in the general setting described below, in the spirit of not missing the forest for the trees by sparing the reader the details of some tedious calculations. This provides the main motivation for this Section.

Lower-case Latin indices belong to the set $\{1,2,3\}$ representing Rectangular Cartesian spatial coordinates, and $t$ is time. Let upper-case Latin indices belong to the set $\{1,2,3, \cdots,N \}$, indexing the components of the $N \times 1$ array of primal variables, $U$, with, possibly, a conversion to first-order form as necessary. We consider the system of equations
\begin{subequations}
    \label{eq:gov_cont_mech}
    \begin{align}
       \calC_{\Gamma I} \p_t U_I + \p_j \calF_{\Gamma j}(U) + G_\Gamma(U,x,t) & = 0 \ \mbox{in} \ \Omega \times (0,T),  \qquad \Gamma = 1, \ldots, N^* \label{eq:gov}\\
        \calC_{\Gamma I} U_I (x, 0) & = \calC_{\Gamma I} U_I^{(0)} (x) \ \mbox{specified on } \Omega \ \mbox{(initial conditions)}\\
        (\calF_{\Gamma j}(U)n_j)|_{(x,t)} & = (B_{\Gamma j} n_j)|_{(x,t)}  \ \mbox{specified on } \p\Omega_\Gamma \ \mbox{(boundary conditions)} \label{eq:gov_bc},
    \end{align}
\end{subequations}
where $\Omega$ is a fixed domain in $\R^3$ with boundary $\p \Omega \supset  \bigcup_\Gamma \p \Omega_\Gamma$, upper-case Greek indices index the number of equations involved, after conversion to first-order form when needed. Here, $\calC$ is an $N^* \times N$ matrix, $\calF, G$ are given functions of their argument, and $U^{(0)}, B$ are specified functions. 

It can be shown that nonlinear elastostatics, elastodynamics, (in)compressible Euler and Navier Stokes can all be written in this form. In this work, we will explicitly consider the cases of classical rate-dependent and rate-independent single crystal plasticity, the latter furnishing a concrete setting for considering inequality constraints, converted to equalities by the addition of slack variables.  As an example, consider the equations of nonlinear elastostatics given by 
\begin{subequations}
    \label{eq:elastostat}
    \begin{align}
        \p_j \hat{P}_{ij}(F)  & = 0 \ \mbox{in} \ \Omega \label{eq:equil}\\
        \p_j y_i - F_{ij} & = 0  \ \mbox{in} \ \Omega \label{eq:y-F}\\
        \hat{P}_{ij}n_j  = p_i  \ \mbox{on} \ \p\Omega_p \qquad & ; \qquad y_i = y^{(b)}_i  \ \mbox{on} \ \p\Omega_y \label{eq:bc_stat}
    \end{align}
\end{subequations}
where $\hat{P}$ is the First Piola-Kirchhoff stress response function. Let $(y_1, y_2, y_3)$ form the first three components of the array $U$.  The conversion to first-order form (so that the Lagrangian $\scl_H$ that appears subsequently in \eqref{eq:defs} contains no derivatives in the primal variables) requires the addition of nine more primal variables $F$ \eqref{eq:y-F}. These additional nine relations can be written in the form
\begin{equation}\label{eq:first-order}
 \mathcal{A}_{\Gamma I j} \p_j U_I - \mathcal{B}_{\Gamma I} U_I = 0, \qquad \Gamma = 4,\ldots 12,   
\end{equation}
where $\mathcal{A}, \mathcal{B}$ are constant matrices (with $\mathcal{B}$ diagonal in many cases) that define the augmentation of the primal list from $(y)$ to $(y, F)$, and define the augmenting primal variables as, in general, linear combinations of the partial derivatives of components of $U$. The equation set \eqref{eq:first-order} can be expressed in  the form \eqref{eq:gov}, and we note, for the convenience of the reader, that the arrays $\mathcal{B}$ and $B$ are not the same.

Boundary conditions are best considered on a specific case-by-case basis. It is shown in Appendix \ref{app:A} how Dirichlet boundary conditions can be accommodated within the setup \eqref{eq:gov_cont_mech}.

Define the pre-dual functional by forming the scalar products of \eqref{eq:gov} with the dual fields $D$, integrating by parts, substituting the prescribed initial and boundary conditions (ignoring, for now, space-time boundary contributions that are not specified)  \textit{and adding a potential $H$ as shown}:
\begin{equation}
    \label{eq:predual}
    \begin{aligned}
        \widehat{S}_H[U,D] & = \int_\Omega \int_0^T \left( - \calC_{\Gamma I} \p_j U_I \p_t D_\Gamma - \calF_{\Gamma j}\big|_U \p_j D_\Gamma + G_\Gamma\big|_{(U,x,t)} D_\Gamma + H(U,x,t) \right) \,dx dt  \\
        & \qquad - \int_\Omega \calC_{\Gamma I} {U}_I^{(0)}(x) D(x, 0) \, dx + \sum_\Gamma \int_{\p \Omega_\Gamma} \int_0^T  B_{\Gamma j} \, D_\Gamma \, n_j \, da dt,
    \end{aligned}    
\end{equation}
(where the arguments $(x,t)$ are suppressed except to display the  explicit dependence of $G, H$ and in the initial condition).

Define
\begin{equation}
    \label{eq:defs}
    \begin{aligned}
        \dee & := (\p_t D, \nabla D, D)\\
        \scl_H(U,\dee,x,t) & := - \calC_{\Gamma I} U_I \p_t D_\Gamma - \calF_{\Gamma j}\big|_U \p_j D_\Gamma + G_\Gamma(U,x,t) D_\Gamma + H(U,x,t)
    \end{aligned}
\end{equation}
and \emph{require the choice} of the potential $H$ to be such that it facilitates the existence of a function
\[
U = U^{(H)}(\dee, x, t)
\]
which satisfies
\begin{equation}\label{eq:dLdU}
\frac{\p \scl_H}{\p U} \left(U^{(H)}(\dee, x,t),\dee,x,t \right) = 0 \qquad \forall \ (\dee,x,t).
\end{equation}
When such a \textit{dual-to-primal} (DtP) `change of variables' mapping, $U^{(H)}$, exists, defining the \textit{dual} functional as
\begin{equation*}
        \begin{aligned}
        & S_H[D] := \widehat{S}_H \left[ U^{(H)}, D \right] \\
        &  \ = \int_\Omega \int_0^T \scl_H\left(U^{(H)}(\dee,x,t), \dee, x, t \right) \, dxdt - \int_\Omega \calC_{\Gamma I} U_I^{(0)}(x) D_\Gamma(x, 0) \, dx + \sum_\Gamma \int_{\p \Omega_\Gamma} \int_0^T  B_{\Gamma j} \, D_\Gamma \, n_j \, da dt,\\
        & \mbox{with $D$ specified (arbitrarily) on parts of the space-time domain boundary complementary to those }\\
        & \mbox{that appear explicitly above,}
    \end{aligned}    
\end{equation*}
and noting \eqref{eq:dLdU}, the first variation of $S_H$ (about a state $(x,t) \mapsto D(x,t)$ in the direction $\delta D$, the latter constrained to vanish on parts of the boundary where $D$ is specified), is given by
\begin{equation*}
\begin{aligned}
        \delta S_H\bigg|_{\delta D} [D] & = \bigintsss_\Omega \bigintsss_0^T \frac{ \p \scl_H}{\p \dee} \left(U^{(H)}(\dee,x,t), \dee, x,t \right) \cdot \delta \dee \, dx dt - \int_\Omega \calC_{\Gamma I} {U}_I^{(0)}(x) \delta D_\Gamma(x, 0) \, dx \\
    & \quad + \sum_\Gamma \int_{\p \Omega_\Gamma} \int_0^T  B_{\Gamma j} \, \delta D_\Gamma \, n_j \, da dt.
\end{aligned}
\end{equation*}

Noting, now, that $\scl_H$ is necessarily affine in $\dee$, its second argument, it can be checked that \textit{the Euler-Lagrange (E-L) equations and natural boundary conditions of the dual functional $S_H$ are exactly the system \eqref{eq:gov} with $U$ substituted by $U^{(H)}(\dee|_{(x,\cdot)},x,\cdot)$}; the first variation is explicitly given as
\begin{equation*}
    \begin{aligned}
        & \delta S_H\Big|_{\delta D} [D]  = \\
        & \int_\Omega \int_0^T \left( \p_t\left( \calC_{\Gamma I} U^{(H)}_I \Big|_{\left (\dee|_{(x,t)},x,t \right)} \right)  + \p_j \calF_{\Gamma j} \left( U^{(H)}\Big|_{\left (\dee|_{(x,t)},x,t \right)}\right) + G_\Gamma \left( U^{(H)}\Big|_{\left (\dee|_{(x,t)},x,t \right)}, x , t \right)\right) \delta D_\Gamma (x,t) \, dx dt\\
        & \quad + \sum_\Gamma \int_{\Omega_\Gamma} \int_0^T \left( B_{\Gamma j} (x,t) - \calF_{\Gamma j} \left( U^{(H)}\Big|_{\left (\dee|_{(x,t)},x,t \right)} \right) \right) n_j (x,t) \delta D_\Gamma(x,t) \, dadt \\
        & \quad + \int_\Omega \calC_{\Gamma I} \left( U^{(H)}_I \Big|_{\left(\dee|_{(x,0)},x,0 \right)} - U^{(0)}_I(x) \right) \delta D (x, 0) \, dx.
    \end{aligned}
\end{equation*}

It is this simple idea that we exploit to develop variational principles for a class of models from continuum mechanics.

It is an important consistency check of our scheme that considering the potential $H$ of the form
\begin{equation}\label{eq:gen_H}
H(U,x,t) = \frac{1}{2}\,  a_U \,\left| U - \bar{U}(x,t) \right|^2 + \frac{1}{p} \, b_U \, \left| U - \bar{U}(x,t) \right|^p,
\end{equation}
where $a_U, b_U$ are positive constants, typically large, with $p >2$ tailored to the nonlinearities present in the functions $\calF, G$, and for $(x,t) \mapsto \bar{U}(x,t)$ an arbitrarily specified function,
\begin{equation} \label{eq:dLdU_powerlaw}
\frac{\p \scl}{\p U_I} = - \calC_{\Gamma I} \p_t D_\Gamma - \frac{ \p \calF_{\Gamma j}}{\p U_I} \p_j D_\Gamma + \frac{\p G_\Gamma}{\p U_I} D_\Gamma + \left( a_U + b_U \, \left|U - \bar{U} \right|^{p-2} \right) \left(U_I - \bar{U}_I \right) = 0
\end{equation} 
is solved, 
\[
\mbox{for} \qquad  \dee(x,t) := (\p_t D(x,t), \nabla D (x,t), D(x,t)) = (0,0,0), \qquad \mbox{by} \qquad U^{(H)}(\dee(x,t),x,t) = \bar{U}(x,t).
\]
\textit{If we now choose} $\bar{U}$ \textit{as a solution to the primal problem} \eqref{eq:gov_cont_mech}, \textit{then a (smooth) solution exists to the E-L equations of the dual problem} given by $(x,t) \mapsto D(x,t) = 0$. This is an existence result for our dual problem. As well, it shows that all solutions to the primal problem can be recovered by the dual scheme by a family of appropriately designed dual problems. 

Of course, it is the goal of our strategy to design and use specific $H$'s, without the knowledge of exact solutions to the primal problem, as a selection criterion to recover special sets of (possibly unstable) solutions of the primal problem in a `stable' manner by solving the dual problem. We note that there are examples in continuum mechanics, e.g.~nonconvex elastostatics in 1-d, where an unstable solution (or critical point) of a primal energy functional is actually the limit of an energy minimizing sequence, which is then recovered as a minima of a relaxed primal problem.

Another important point to note is that the dual E-L equations corresponding to primal \textit{initial}-(boundary)-value problems contain second order time derivatives in the dual variables, after conversion of the primal system to first-order form; this can be understood by considering the form of the DtP mapping \eqref{eq:dLdU_powerlaw} and the primal system \eqref{eq:gov}. This generally requires two `boundary' conditions in the time-like direction on such variables, when at most, only one is available from the primal problem. This raises the question of how the second condition ought to be specified and what effect it has on the recovery of the correct primal solution, especially when the primal system has a unique solution as an initial-value problem. It turns out that a final time boundary condition can be arbitrarily specified on the dual variables and this does not have an effect on the recovery of correct primal solutions as the DtP mapping, for standard initial value problems, necessarily depends on $\p_t D$, see e.g.~\eqref{eq:dLdU_powerlaw}, and specifying $D$ at the final time leaves the time derivative free to adjust to the demands of achieving the required primal solution through the DtP mapping. This fact has been discussed and demonstrated in specific contexts in \cite[Sec.~7]{action_2} and \cite{KA1}.

\section{Local degenerate ellipticity of the dual formulation of continuum mechanics}\label{sec:deg_ellipt}
For this section, let Greek lower-case indices belong to the set $\{0,1,2,3\}$ representing Rectangular Cartesian space-time coordinates $x^\alpha, \alpha = 0, 1,2,3$; $0$ represents the time coordinate when the PDE is time-dependent. Let upper-case Latin indices belong to the set $\{1,2,3, \cdots,N \}$, indexing the components of the $N \times 1$ array of primal variables, $U$, with, possibly, a conversion to first-order form as necessary. Now consider the system of primal PDE
\begin{equation}\label{eq:conserv_source}
    \p_\alpha (\calF_{\Gamma\alpha}(U)) + G_\Gamma(U, x) = 0, \qquad \Gamma = 1, \ldots, N^*
\end{equation}
where upper-case Greek indices index the number of equations involved, after conversion to first-order form when needed.  

We assume that the functions $U \mapsto \frac{\p^2 \calF_\Gamma}{\p U_P \p U_R}(U)$ and $U \mapsto \frac{\p^2 G_\Gamma}{\p U_P \p U_R}(U)$ are bounded functions on their domains.

Let $D$ be the $N^* \times 1$ array of dual fields and, as earlier, let us consider a shifted quadratic for the potential $H$, characterized by a diagonal matrix $[a_{kj}]$ with constant positive diagonal entries so that the Lagrangian takes the form (with $H$ chosen as a quadratic form for simplicity of presentation)
\begin{equation*}
    \scl(U,D,\nabla D, \bar{U}) := - \calF_{\Gamma \alpha}(U) \p_\alpha D_\Gamma + D_\Gamma G_\Gamma(U) + \half (U_k - \bar{U}_k) a_{kj} (U_j - \bar{U}_j).
\end{equation*}
Then the corresponding DtP mapping, obtained by `solving $\frac{\p \scl}{\p U} = 0$ for $U$ in terms of $(\nabla D, D, \bar{U})$,' is given by the implicit equation
\begin{equation}\label{eq:deg_ell_map}
    U_J^{(Q)}(\nabla D, D, \bar{U}) = \bar{U}_J + (a^{-1})_{JK} \left( \frac{\p \calF_{\Gamma \alpha}}{\p U_K} \bigg|_{U^{(Q)}(\nabla D, D, \bar{U})} \p_\alpha D_\Gamma - D_\Gamma \frac{\p G_\Gamma}{\p U_K}\bigg|_{U^{(Q)}(\nabla D, D, \bar{U})} \right).
\end{equation}
It is a fundamental property of the dual scheme that the dual E-L equation is then given by
\begin{equation}\label{eq:dual_EL}
    \p_\alpha \Big( \calF_{\Gamma \alpha} \big( U(\nabla D, D, \bar{U})) \Big) + G_\Gamma ( U(\nabla D, D, \bar{U}) ) = 0
\end{equation}
(where we have dropped the superscript $^{(Q)}$ for notational convenience), whose ellipticity is governed by the term
\[
\mathbb{A}_{\Gamma \alpha \Pi \mu}(\nabla D, D, \bar{U}) := \frac{\p \calF_{\Gamma \alpha}}{\p U_P}\bigg|_{U(\nabla D, D, \bar{U})} \, \frac{\p U_P}{\p (\nabla D)_{\Pi \mu}}\bigg|_{U(\nabla D, D, \bar{U})}.
\]
From \eqref{eq:deg_ell_map} we have
\begin{equation*}
    \begin{aligned}
         a^{-1}_{PR} \left(  \delta_{\Gamma \Pi} \delta_{\mu \alpha} \frac{ \p \calF_{\Gamma \alpha}}{\p U_R} \, + \, \p_\alpha D_\Gamma \frac{ \p^2 \calF_{\Gamma \alpha}}{\p U_R \p U_S} \frac{\p U_S}{\p (\nabla D)_{\Pi \mu}}  - D_\Gamma \frac{ \p^2 G_\Gamma}{\p U_R \p U_S} \frac{\p U_S}{\p (\nabla D)_{\Pi \mu}} \right)     & =  \frac{\p U_P}{\p (\nabla D)_{\Pi \mu}} \\
        \Longrightarrow \left( \delta_{PS} - a^{-1}_{PR} \p_\alpha D_\Gamma \frac{ \p^2 \calF_{\Gamma \alpha}}{\p U_R \p U_S} + a^{-1}_{PR} D_\Gamma \frac{ \p^2 G_\Gamma}{\p U_R \p U_S} \right) \frac{\p U_S}{\p (\nabla D)_{\Pi \mu}} & = a^{-1}_{PR} \frac{\p \calF_{\Pi \mu}}{\p U_R},
    \end{aligned}
\end{equation*}
and so
\begin{equation*}
    \mathbb{A}_{\Gamma \alpha \Pi \mu}(0,0,\bar{U}) = \frac{\p \calF_{\Gamma \alpha}}{\p U_P}\bigg|_{\bar{U}} \, a^{-1}_{PR} \, \frac{\p \calF_{\Pi \mu}}{\p U_R}\bigg|_{\bar{U}},
\end{equation*}
which is \textit{positive semi-definite} on the space of $N^* \times 3$ (or $N^* \times 4$) matrices. This establishes the degenerate ellipticity of the dual system at the state $x \mapsto D(x) = 0$.

To examine the ellipticity-related properties of the system in a bounded neighborhood, say $\mathcal{N}$, of $(D = 0, \nabla D = 0) \in \mathbb{R}^{N^*} \times \mathbb{R}^{N^* \times \bar{\alpha}}, \bar{\alpha} = 3,4$, we define
\[
M_{PS} := \delta_{PS} - a^{-1}_{PR} \left( \p_\alpha D_\Gamma \frac{\p^2 \calF_{\Gamma \alpha}}{\p U_R \p U_S} - D_\Gamma  \frac{\p^2 G_\Gamma}{\p U_R \p U_S} \right),
\]
and note that
\[
\frac{\p U_P}{\p (\nabla D)_{\Pi \mu} } = M^{-1}_{PQ} a^{-1}_{QR} \frac{\p \calF_{\Pi \mu}}{\p U_R},
\]
where $M^{-1}$ exists and is positive definite by the boundedness of $\mathcal{N}$ and the second derivatives of the functions $\calF$ and $G_\Gamma$, along with an appropriately large choice of the elements of the diagonal matrix $[a_{ij}]$ (in case the second-derivatives are not bounded in some regions of the domain of primal variables we assume that the functions are such that the positive-definiteness of $M$ is maintained. Alternatively, the choice of $H$ can be enhanced (as, e.g.~in \eqref{eq:gen_H}) to dominate the growth of the second derivatives, catering to the specifics of the second-derivative functions in any particular problem).

The degenerate ellipticity or `convexity' of the system \eqref{eq:conserv_source} in the neighborhood $\mathcal{N}$ is now defined as the positive semi-definiteness of the matrix $\mathbb{A}$  on the space $\mathbb{R}^{N^* \times \bar{\alpha}}$ of matrices, and this in turn is governed by the matrix
\[
\mathbb{A}^{(sym)}_{\Gamma \alpha \Pi \mu}\bigg|_{(\nabla D, D, \bar{U})} = \frac{\p \calF_{\Gamma \alpha}}{\p U_P}\bigg|_{U(\nabla D, D, \bar{U})} \frac{1}{2} \left( M^{-1}_{PQ}\bigg|_{U(\nabla D, D, \bar{U})}  a^{-1}_{QR} + M^{-1}_{RQ}\bigg|_{U(\nabla D, D, \bar{U})}  a^{-1}_{QP}\right)\frac{\p \calF_{\Pi \mu}}{\p U_R}\bigg|_{U(\nabla D, D, \bar{U})}.
\]
By the positive definiteness of the matrix $[M_{PS}]$ in the neighborhood $\mathcal{N}$, it follows that
\[
\xi^{\Gamma \alpha} \ \mathbb{A}^{(sym)}_{\Gamma \alpha \Pi \mu}\bigg|_{(\nabla D, D, \bar{U})} \xi^{\Pi \mu}  \geq 0 \qquad \forall \qquad (D, \nabla D) \in \mathcal{N}, \quad \xi \in \R^{N^* \times \bar{\alpha}}
\]
which establishes a `local' degenerate ellipticity of the system \eqref{eq:conserv_source}. We note that degenerate ellipticity is stronger than the Legendre-Hadamard condition given by the requirement of positive semi-definiteness of $\mathbb{A}$ on the space of tensor products from $\mathbb{R}^{N^*} \otimes \mathbb{R}^{\bar{\alpha}}$, and not directly comparable to the strong-ellipticity condition, since it is weaker than the latter when restricted to the space $\mathbb{R}^{N^*} \otimes \mathbb{R}^{\bar{\alpha}}$ but simultaneously requiring semi-definiteness on the larger space of $\mathbb{R}^{N^* \times \bar{\alpha}}$. Also of note is that degenerate ellipticity does not preclude the failure of ellipticity characterized by the condition $\det [\mathbb{A}_{\Gamma \alpha \Pi \mu} n_\alpha n_\mu] \neq 0$ for all unit direction $n \in \R^ {\bar{\alpha}}$, $\bar{\alpha} = 3$ or $4$, thus allowing for weak (gradient) discontinuities of weak solutions $x \mapsto D(x)$ of \eqref{eq:dual_EL} (or at least its linearized counterpart), a feature that is important for recovering discontinuous solutions of the primal problem (e.g.~inviscid Burgers) expressed as combinations of derivatives of the dual fields through the DtP mapping as, e.g., demonstrated in the context of the linear transport equation in \cite{KA1}.

If a solution of the primal system is close to the base state $\bar{U}$, then it seems natural to  expect, due to this local degenerate ellipticity, that such a solution can be obtained in a `stable' manner by the dual formulation \textit{designed by the choice of the auxiliary potential $H$ as a shifted quadratic (or `power law') about the base state $\bar{U}$}, for instance by an iterative scheme starting from a guess $(D = 0, U = \bar{U})$.

Our experience \cite{KA1, sga, KA2, a_arora_thesis} shows that this observation is of great practical relevance in using the dual scheme, and we consistently exploit it in all our computational approximations.

To make contact with the parlance of the classical `rate problems' of Hill \cite{hill1957uniqueness, hill1956new, hill1979aspects},  degenerate ellipticity here corresponds to the absence of negative `energy' modes of the linearized, or `incremental/rate,' dual problem at dual states whose corresponding primal state, obtained via the DtP mapping, may well entail a loss of positive-semi-definiteness of the physical incremental moduli on the space of dyads $a \otimes n \ ( a, n \in \R^{3})$ in the primal rate problem under quasi-static conditions. 

Furthermore, by a theorem of Ball \cite{ball1976convexity} and in the context of nonlinear hyperelasticity as the primal problem, quasiconvexity implies the Legendre-Hadamard condition (for the primal problem) so that it is possible that the dual problem remains degenerate elliptic/convex, even when the primal problem is not quasiconvex.

\section{A variational principle for rate-dependent, dynamic, single crystal plasticity}\label{sec:rd_sc}
We follow the scheme described in Sec.~\ref{sec:dual_cont_mech} to develop the required variational principle. The specifics of rate-dependent single crystal plasticity theory can be found in the expositions  of \cite{hutchinson1976bounds, asaro1983micromechanics}. 

Let $\Omega \subset \R^3$ be a given, fixed reference configuration with all spatial derivatives below being w.r.t rectangular Cartesian coordinates parametrizing this reference, and partial derivatives w.r.t time, $t$, representing  material time derivatives, also alternatively written with a superposed dot. The interval $[0,T]$ is fixed, but chosen arbitrarily. Lowercase Greek (super)subscripts refer to numbering of slip systems. We consider the following set of equations on $\Omega$:
\begin{equation}\label{eq:primal_rd}
    \begin{aligned}
        \rho_0 \dot{v} - \p_j \calN_{ij}(F,P) & = 0 \\
        \dot{P}_{ij} - \sum_\alpha \Big( r^\alpha(F,P,g) \, m^\alpha_i  n^\alpha_k \Big) P_{kj} & = 0 \\
        \dot{g}_\alpha - h_{\alpha \beta} (g) \, r^\beta(F, P, g) & = 0 \\
        \dot{y}_i - v_i & = 0 \\
        \p_j y_i - F_{ij} & = 0,
    \end{aligned}
\end{equation}
with the boundary conditions
\begin{equation}\label{eq:primal_rd_bc}
    \calN_{ij}(F,P)\big|_{(x,t)} n_j\big|_x = \bar{t}_i(x,t), \  x \in \p \Omega_{\bar{t}}; \qquad y_i(x,t) = y^{(b)}_i(x,t), \ x \in \p \Omega_y,
\end{equation}
and initial conditions
\begin{equation}\label{eq:eq:primal_rd_ic}
y_i(x,0) = y^{(0)}_i(x), \qquad v_i(x,0) = v^{(0)}_i(x), \qquad P_{ij} (x,0) = P^{(0)}_{ij} (x), \qquad g^{\alpha}(x,0) = g^{\alpha(0)}(x), \qquad x \in \Omega.
\end{equation}

In the above, $\rho_0$ is a given mass density field on the reference configuration, $\calN$ is the response function for the first Piola-Kirchhoff stress w.r.t.~the reference configuration, $y,v, F$ are the position, velocity, and deformation gradient fields, respectively, $P$ is the plastic distortion tensor, $r^\alpha$ are response functions for the slip system rates (e.g., the power law \cite{hutchinson1976bounds} or the Perzyna overstress model \cite{perzyna1966fundamental}), $(m^\alpha, n^\alpha)$ are the elastically unstretched slip direction and slip normal vectors, $g^\alpha$ are the strengths, and $h_{\alpha \beta}$ are the hardening matrix response functions \cite{hill1966generalized}. All quantities indexed by $\alpha$ refer to an object corresponding to the $\alpha^{th}$ slip system. The functions $\bar{t},y^{(b)}, y^{(0)}, v^{(0)}, P^{(0)}, g^{\alpha(0)}$ are prescribed.

Now define the array of primal fields
\[
U = (y, v, F, P, g),
\]
the dual fields
\[
D = (\xi, \gamma, \Phi, \Pi, \Gamma),
\]
and assume the potential $H$ to be of the form
\begin{equation*}
    \begin{aligned}
        & H(y,v,F,P, g, x, t) = \\
        & \quad \frac{1}{2} \left( a_y \Big|y - \bar{y} |_{(x,t)} \Big|^2 + a_v \Big|v - \bar{v}|_{(x,t)} \Big|^2  + a_F \Big|F - \bar{F} |_{(x,t)} \Big|^2 + a_P \Big|P - \bar{P}|_{(x,t)} \Big|^2 + a_g \Big|g - \bar{g}|_{(x,t)} \Big|^2 \right) \\
        &  + \frac{1}{p} \left(  b_F \big|F - \bar{F}|_{(x,t)} \Big|^p + b_P \Big|P - \bar{P}|_{(x,t)} \Big|^p + b_g \Big|g - \bar{g}|_{(x,t)} \Big|^p \right),
    \end{aligned}
\end{equation*}
for $p > 2$ as needed. 

Here, the \emph{base states}, the collection of space-time fields with overhead bars, are arbitrarily specified, with their closeness to an actual solution of the primal problem resulting in a better design of the variational principle. 

The introduction of the power $p$ is simply to ensure strict convexity of the ensuing Lagrangian $\scl_H$ in the primal variables $U$ for each fixed set of values of $(\dee, x, t)$, which in turn is closely dictated by the nonlinearities of the primal problem. The non-negative real-valued constants $a_{(\cdot)}, b_{(\cdot)}$ are chosen arbitrarily, typically large, when non-zero, to facilitate the strict convexity of $\scl_H$ that appears below. Clearly, there is a great deal of freedom in making the choices of $H$; e.g., the power $p$ in the specific choice above does not even have to be the same on each of the terms.

We now define the pre-dual functional
\begin{equation}\label{eq:S_hat_rd}
    \begin{aligned}
        \widehat{S}[U, D] & =  \int_\Omega \int_0^T \scl_H (U, \dee,x, t) \, dx dt + \mbox{inital and boundary contributions} \\
        & := \ \int_\Omega \int_0^T  - \rho_0 v_i \p_t \gamma_i + \calN_{ij}\big|_{(F,P)} \p_j \gamma_i  - P_{ij} \p_t \Pi_{ij} - \Pi_{ij} \sum_\alpha r^\alpha \big|_{(F,P,g)}  \, m^\alpha_i  n^\alpha_k P_{kj}  \  dx dt \\  
        & \quad - \int_\Omega \rho_0 v^{(0)}_i(x) \gamma_i(x,0)  + P^{(0)}_{ij}(x) \Pi_{ij}(x,0) \, dx - \int_{\p \Omega_{\bar{t}}} \int_0^T \bar{t}_i \gamma_i \, da dt \\
        & \quad - \int_\Omega \int_0^T  g_\alpha \p_t \Gamma^\alpha  + \Gamma^\alpha h_{\alpha \beta}\big|_g r^\beta \big|_{(F,P,g)} + y_i \p_t \xi_i + \xi_i v_i + y_i \p_j \Phi_{ij} + \Phi_{ij} F_{ij}  \, dx dt \\
        & \quad - \int_\Omega g^{\alpha (0)}(x) \Gamma^\alpha(x,0) + y_i^{(0)}(x) \xi_i(x,0) \, dx + \int_{\p \Omega_y} \int_0^T y_i^{(b)} \Phi_{ij} n_j \, da dt\\
        & \quad + \int_\Omega \int_0^T H(y,v,F,P, g, x, t) \, dxdt ,
    \end{aligned}
\end{equation}
where the array $\dee$ is defined as
\[
\dee = (\p_t \gamma, \nabla \gamma, \p_t \Pi, \Pi, \p_t \Gamma, \Gamma, \p_t \xi, \xi, div \, \Phi,  \Phi).
\]
In order to define the function $U^{(H)}$ we need to consider the $(x,t)$-pointwise equations for $U^{(H)}$ for the given values $(\dee(x,t), \bar{U}(x,t))$ (we will drop the superscript $^{(H)}$ for notational convenience):
\begin{equation*}
    \begin{aligned}
       & \frac{\p \scl_H}{\p y_i} : \qquad - \p_t \xi_i - \p_j \Phi_{ij} + a_y \left( y_i - \bar{y}_i \right) = 0\\
       & \frac{\p \scl_H}{\p v_i}  : \qquad - \rho_0 \p_t \gamma_i - \xi_i + a_v \left( v_i - \bar{v}_i \right) = 0\\
       & \frac{\p \scl_H}{\p F_{ij}} : \qquad \p_l \gamma_k \frac{\p \calN_{kl}}{\p F_{ij}}\bigg|_{(F,P)} - \Pi_{rs} \sum_\alpha \frac{\p r^\alpha}{\p F_{ij}}\bigg|_{(F,P,g)} m^\alpha_r n^\alpha_k P_{ks}  \\
       & \qquad \qquad \quad - \Gamma^\alpha h_{\alpha \beta}|_g \frac{\p r^\beta}{\p F_{ij}}\bigg|_{(F,P,g)} - \Phi_{ij} + \left(a_F + b_F \left| F - \bar{F} \right|^{p-2}\right) \left( F_{ij} - \bar{F}_{ij} \right) = 0\\
       & \frac{\p \scl_H}{\p P_{rs}} : \qquad \p_l \gamma_k \frac{\p \calN_{kl}}{\p P_{rs}}\bigg|_{(F,P)} - \p_t \Pi_{rs} - \Pi_{ij} \sum_\alpha \frac{\p r^\alpha}{\p P_{rs}}\bigg|_{(F,P,g)} m^\alpha_i n^\alpha_k P_{kj} \\
       & \qquad \qquad \quad  - \Pi_{is} \sum_\alpha r^\alpha|_{(F,P,g)} m^\alpha_i n^\alpha_r - \Gamma^\alpha h_{\alpha \beta}|_g \frac{\p r^\beta}{\p P_{rs}}\bigg|_{(F,P,g)} + \left(a_P + b_P \left| P - \bar{P} \right|^{p-2}\right) \left( P_{ij} - \bar{P}_{ij} \right) = 0\\
       & \frac{\p \scl_H}{\p g^\alpha} : \qquad - \Pi_{ij} \sum_\kappa \frac{\p r^\kappa}{\p g^\alpha}\bigg|_{(F,P,g)} m^\kappa_i n^\kappa_k P_{kj} - \p_t \Gamma^\alpha - \Gamma^\kappa h_{\kappa \beta}|_g \frac{\p r^\beta}{\p g^\alpha}\bigg|_{(F,P,g)} + \left(a_g + b_g \left| g - \bar{g} \right|^{p-2}\right) \left( g^\alpha - \bar{g}^\alpha \right) = 0.
    \end{aligned}
\end{equation*}
By making suitable choices for the various constants appearing in $H$, the expectation is that $\scl_H$ can be made strictly convex in $U^{(H)}$ so that a (unique) solution for $U^{(H)}$ exists and can be solved by standard techniques without difficulty at (almost) every point of the domain.

We now define a \emph{dual} functional for dynamic rate-dependent, single crystal plasticity as
\[
S[D] = \widehat{S} \left[ U^{(H)}, D\right] 
\]
interpreted as replacing all occurrences of $U$ in the right-hand-side of \eqref{eq:S_hat_rd} by $U^{(H)}\left(\dee, \bar{U}(x,t) \right)$ subject to the following essential `boundary conditions' on parts of the space-time domain boundary given by $\calB := (\p \Omega   \times (0,T)) \cup (\Omega \times \{0,T\})$:
\begin{equation*}
    \begin{aligned}
       &  \mbox{ (arbitrarily) specified $\gamma$ on $\calB \backslash ( \ (\p \Omega_{\bar{t}} \times (0,T)) \cup (\Omega \times \{0\}) \ )$  and $\Phi$ on $\calB \backslash (\p \Omega_y \times (0,T)) $ and} \\
        & \mbox{ (arbitrarily) specified $\Pi, \Gamma, \xi$ on $\calB \backslash (\Omega \times \{0\}) $}.
    \end{aligned}
\end{equation*}
The Euler-Lagrange equations and the natural boundary conditions of $S$ are the equations \eqref{eq:primal_rd}-\eqref{eq:primal_rd_bc}-\eqref{eq:eq:primal_rd_ic}, with the replacement $U \rightarrow U^{(H)}\left(\dee, \bar{U}(x,t) \right)$.

\section{A variational principle for rate-independent, quasi-static, single crystal plasticity}\label{sec:ri_sc}
The specifics of rate-independent single crystal plasticity theory may be found in the expositions  \cite{havner1992finite, bassani1993plastic}.

\textit{In this section, the summation convention is not used on lower-case Greek indices} which index the slip systems. We consider the following set of equations on a fixed reference $\Omega \subset \R^3$:
\begin{equation}\label{eq:primal_ri}
    \begin{aligned}
        \p_j \calN_{ij}(F,P) & = 0 \\
        \dot{P}_{ij} - \sum_\alpha \Big( r^\alpha \, m^\alpha_i  n^\alpha_k \Big) P_{kj} & = 0 \\
        \dot{g}_\alpha - \sum_\beta h_{\alpha \beta} (g) \, r^\beta(F, P, g) & = 0 \\
        \p_j y_i - F_{ij} & = 0 \\
        Y^\alpha (F,P,g) + s_\alpha^2 & = 0 \\
        r^\alpha Y^\alpha & = 0 \\
        r^\alpha - p_\alpha^2 & = 0,
    \end{aligned}
\end{equation}
with the boundary conditions
\begin{equation}\label{eq:primal_ri_bc}
    \calN_{ij}(F,P)\big|_{(x,t)} n_j\big|_x = \bar{t}_i(x,t), \  x \in \p \Omega_{\bar{t}}; \qquad y_i(x,t) = y^{(b)}_i(x,t), \ x \in \p \Omega_y,
\end{equation}
and initial conditions
\begin{equation}\label{eq:eq:primal_ri_ic}
P_{ij} (x,0) = P^{(0)}_{ij} (x), \qquad g^{\alpha}(x,0) = g^{\alpha(0)}(x), \qquad x \in \Omega.
\end{equation}
In the above, $\calN$ is the response function for the first Piola-Kirchhoff stress w.r.t. the reference configuration, $y, F$ are the position and deformation gradient fields, respectively, $P$ is the plastic distortion tensor, $r^\alpha$ is a slip-rate, $(m^\alpha, n^\alpha)$ are the unstretched slip direction and slip normal vectors, $g^\alpha$ is a strength, $h_{\alpha \beta}$ are the hardening matrix response functions, $Y^\alpha$ is an yield response function (the canonical example being $Y^\alpha = \tau^\alpha - g^\alpha$, where $\tau^\alpha$ is the resolved shear stress on the slip system $\alpha$ given by $\tau^\alpha = (F^e m^\alpha)_i T_{ij} (F^{e-T} n^\alpha)_j$ where $F^e := FP^{-1}$ is the elastic distortion, and $T = (\det F)^{-1} \calN F^T$ is the Cauchy stress tensor), and $s_\alpha, p_\alpha$ are slack variables. All quantities indexed by $\alpha$ refer to an object corresponding to the $\alpha^{th}$ slip system. The slack variables enable the imposition of the inequalities
\[
Y^\alpha \leq 0; \qquad r^\alpha \geq 0.
\]
 The functions $\bar{t},y^{(b)}, P^{(0)}, g^{\alpha(0)}$ are prescribed.

Now define the array of primal fields
\[
U = (y, F, P, g, r, s, p),
\]
the dual fields
\[
D = (\gamma, \Phi, \Pi, \Gamma, \mu, \rho, \nu),
\]
and assume the potential $H$ to be of the form
\begin{equation*}
    \begin{aligned}
        & H(y,v,F,P, g, r,s, p, x, t) = \\
        & \quad \frac{1}{2} \left( a_y \Big|y - \bar{y} |_{(x,t)} \Big|^2 + a_F \Big|F - \bar{F} |_{(x,t)} \Big|^2 + a_P \Big|P - \bar{P}|_{(x,t)} \Big|^2 + a_g \Big|g - \bar{g}|_{(x,t)} \Big|^2 \right) \\
        &  + \frac{1}{2} \left(  a_r \big|r - \bar{r}|_{(x,t)} \big|^2 + a_s \big|s - \bar{s}|_{(x,t)} \big|^2 + a_p \big|p - \bar{p}|_{(x,t)} \big|^2 \right)\\
        &  + \frac{1}{p} \left(  b_F \Big|F - \bar{F}|_{(x,t)} \Big|^p + b_P \Big|P - \bar{P}|_{(x,t)} \Big|^p + b_g \Big|g - \bar{g}|_{(x,t)} \Big|^p \right),
    \end{aligned}
\end{equation*}
for $p > 2$ as needed, with the same understanding operative for base states and the various constants that appear as in the previous Section \ref{sec:rd_sc}.

We now define the pre-dual functional
\begin{equation}\label{eq:S_hat_ri}
    \begin{aligned}
        \widehat{S}[U, D] & =  \int_\Omega \int_0^T \scl_H (U, \dee,x, t) \, dx dt + \mbox{inital and boundary contributions} \\
        & := \ \int_\Omega \int_0^T  - \calN_{ij}\big|_{(F,P)} \p_j \gamma_i  - P_{ij} \p_t \Pi_{ij} - \Pi_{ij} \sum_\alpha r^\alpha  \, m^\alpha_i  n^\alpha_k P_{kj}  \  dx dt \\  
        & \quad  - \int_\Omega P^{(0)}_{ij}(x) \Pi_{ij}(x,0) \, dx + \int_{\p \Omega_{\bar{t}}} \int_0^T \bar{t}_i \gamma_i \, da dt + \int_{\p \Omega_y} \int_0^T y_i^{(b)} \Phi_{ij} n_j \, da dt\\
        & \quad - \int_\Omega \int_0^T y_i \p_j \Phi_{ij} + \Phi_{ij} F_{ij} + \sum_\alpha g^\alpha \p_t \Gamma^\alpha  + \sum_\alpha \sum_\beta \Gamma^\alpha h_{\alpha \beta}\big|_g r^\beta    \, dx dt \\
        & \quad +  \int_\Omega \int_0^T \sum_\alpha \left( \rho^\alpha Y^\alpha \big|_{(F,P,g)} + \rho^\alpha s_\alpha^2 + r^\alpha Y^\alpha\big|_{(F,P,g)} \mu^\alpha + r^\alpha \nu^\alpha  - \nu^\alpha p_\alpha^2 \right) \, dx dt\\
        & \quad - \int_\Omega g^{\alpha (0)}(x) \Gamma^\alpha(x,0) \, dx \\
        & \quad + \int_\Omega \int_0^T H(y,v,F,P, g, r,s, p, x, t) \, dx dt,
    \end{aligned}
\end{equation}
where the array $\dee$ is defined as
\[
\dee = (\nabla \gamma, \p_t \Pi,  \Pi, div \, \Phi, \Phi, \rho, \mu, \nu, \p_t \Gamma, \Gamma).
\]
In order to define the function $U^{(H)}$ we need to consider the following $(x,t)$-pointwise equations for $U^{(H)}$ for the given values $(\dee(x,t), \bar{U}(x,t))$ (we will drop the superscript $^{(H)}$ for notational convenience):
\begin{equation*}
    \begin{aligned}
       & \frac{\p \scl_H}{\p y_i} : \qquad - \p_j \Phi_{ij} + a_y \left( y_i - \bar{y}_i \right) = 0\\
       & \frac{\p \scl_H}{\p F_{rs}} : \qquad - \p_j \gamma_i \frac{\p \calN_{ij}}{\p F_{rs}}\bigg|_{(F,P)} -  \Phi_{ij}  + \sum_\alpha (\rho^\alpha + r^\alpha \mu^\alpha) \frac{\p Y^\alpha}{\p F_{rs}}\bigg|_{(F,P,g)} + \left(a_F + b_F \left| F - \bar{F} \right|^{p-2}\right) \left( F_{rs} - \bar{F}_{rs} \right) = 0\\
       & \frac{\p \scl_H}{\p P_{rs}} : \qquad - \p_l \gamma_k \frac{\p \calN_{kl}}{\p P_{rs}}\bigg|_{(F,P)} - \p_t \Pi_{rs} - \Pi_{is} \sum_\alpha r^\alpha m^\alpha_i n^\alpha_r + \sum_\alpha (\rho^\alpha + r^\alpha \mu^\alpha) \frac{\p Y^\alpha}{\p P_{rs}}\bigg|_{(F,P,g)} \\
       & \qquad \qquad \quad  + \left(a_P + b_P \left| P - \bar{P} \right|^{p-2}\right) \left( P_{rs} - \bar{P}_{rs} \right) = 0\\
       & \frac{\p \scl_H}{\p g^\mu} : \qquad \sum_\alpha (\rho^\alpha + r^\alpha \mu^\alpha) \frac{\p Y^\alpha}{\p g^\mu}\bigg|_{(F,P,g)} - \p_t \Gamma^\mu - \sum_\alpha \sum_\beta \Gamma^\alpha \frac{\p h_{\alpha \beta}}{\p g^\mu}\bigg|_g r^\beta + \left(a_g + b_g \left| g - \bar{g} \right|^{p-2}\right) \left( g^\mu - \bar{g}^\mu \right) = 0\\
       & \frac{\p \scl_H}{\p r^\alpha} : \qquad - \Pi_{ij} m^\alpha_i n^\alpha_k P_{kj} + Y^\alpha\big|_{(F,P,g)} \mu^\alpha + \nu^\alpha - \sum_\kappa \Gamma^\kappa h_{\kappa \alpha}\big|_g + a_r (r^\alpha - \bar{r}^\alpha) = 0\\
       & \frac{\p \scl_H}{\p s^\alpha} : \qquad 2 s_\alpha \rho^\alpha + a_s (s^\alpha - \bar{s}^\alpha) = 0 \\
       & \frac{\p \scl_H}{\p p^\alpha} : \qquad - 2 p_\alpha \nu^\alpha + a_p (p^\alpha - \bar{p}^\alpha) = 0.
    \end{aligned}
\end{equation*}
Again, by making suitable choices for the various constants appearing in $H$, $\scl_H$ can be made strictly convex in $U^{(H)}$ so that a (unique) solution for $U^{(H)}$ exists and can be solved by standard techniques without difficulty at (almost) every point of the domain.

We now define a \emph{dual} functional for quasi-static, rate-dependent, single crystal plasticity as
\[
S[D] = \widehat{S} \left[ U^{(H)}, D\right] 
\]
interpreted as replacing all occurrences of $U$ in the right-hand-side of \eqref{eq:S_hat_ri} by $U^{(H)}\left(\dee, \bar{U}(x,t) \right)$ subject to the following essential `boundary conditions' on parts of the space-time domain boundary given by $\calB := (\p \Omega   \times (0,T)) \cup (\Omega \times \{0,T\})$:
\begin{equation*}
    \begin{aligned}
       &  \mbox{ (arbitrarily) specified $\gamma$ on $\calB \backslash ( \ (\p \Omega_{\bar{t}} \times (0,T)) \cup (\Omega \times \{0\}) \ )$  and $\Phi$ on $\calB \backslash (\p \Omega_y \times (0,T)) $ and} \\
        & \mbox{ (arbitrarily) specified $\Pi, \Gamma$ on $\calB \backslash (\Omega \times \{0\}) $}.
    \end{aligned}
\end{equation*}
The Euler-Lagrange equations and the natural boundary conditions of $S$ are the equations \eqref{eq:primal_ri}-\eqref{eq:primal_ri_bc}-\eqref{eq:eq:primal_ri_ic}, with the replacement $U \rightarrow U^{(H)}\left(\dee, \bar{U}(x,t) \right)$.

\section{Concluding remarks and outlook}\label{sec:concl}
A formal scheme for developing variational principles for systems of nonlinear partial differential equations arising in continuum mechanics has been proposed.
It is based on the realization that such a system of equations may be viewed as an `invariant' or a `symmetry' of a family of dual variational principles parametrized by a set of scalar potentials of the primal variables, the parametrization acting as the symmetry operation, and the invariant being the Euler-Lagrange equations of any of the variational principles in that family. 

The scheme appears to be best suited for problems which are difficult to solve in the `primal' setting, be it due to lack of existence of solutions as defined by extant strategies, uniqueness, or stability, cf. \cite{a_arora_thesis, sga}, and not meant as a competitor for problems that are solved robustly by existing techniques for the primal problem. It offers the possibility of defining the notion of a very weak solution of the primal problem as the solution to the dual variational problem, which with enough regularity, defines a genuine weak solution of the primal PDE system. This can be useful, as nonlinear PDE systems are generally much harder to solve than a variational minimization/maximization problem. The Euler-Lagrange equations of the dual problem have a certain degenerate ellipticity, and knowledge of `base states' close to  desired solutions can be incorporated in the scheme without approximation; these two features combined together help in obtaining (un)stable solutions of the primal problem in a stable way within the dual formulation - a case study is provided in \cite{sga}. Degenerate ellipticity by itself is not a very strong property (depending on taste, e.g.~when compared to strong ellipticity or strict convexity when physically natural), but does take on significance when the primal PDE system loses ellipticity (along with becoming indefinite) or hyperbolicity.

As a (non-rigorous) sketch of how our scheme may have the potential of achieving the above objective, consider the case of nonlinear hyperelasticity without higher-order regularization. The dominant (and, perhaps, only) strategy available \cite{ball1976convexity} is to declare minimizers of the elastic energy as solutions to the problem. It is well-understood that, more or less, quasiconvexity of the functional (along with some coercivity) is equivalent to the existence of minimizers. As laid out in many works, quasiconvexity is hard to check, 
and it is known to fail for many physical energy functionals 
whose limits of energy minimizing sequences have no status as minimizers of such energy functionals due to the lack of lower semicontinuity.
Juxtaposing the present scheme with this approach, it seeks to define some notion of a solution to the PDE of elastostatics given an elastic stress response function (a system that may be the formal Euler-Lagrange equations of the physical energy functional) thus `severing' the link with looking for minimizers of the physical energy functional and hence its quasiconvexity - and produces a convex variational principle on the dual side, whose critical points and minimizers can nevertheless be sought and, with sufficient regularity in them and the DtP mapping, be deemed as solutions to the PDE of elastostatics - such solutions, of course, need not have a connection to being minimizers of the primal energy functional. Perhaps more importantly, even when the regularity-related steps cannot be carried through, the obtained critical points of the dual problem can be declared as some sort of very weak solutions of the primal PDE system, because of their consistency with the primal problem in the presence of regularity. These steps have been carried out in \cite{sga}.

The formalism has been used in this paper to present variational principles for a class of single crystal plasticity problems which demonstrate its relevance to the theory of generally non-associated, multi-surface plasticity. A particular spin-off of the approach is a potentially robust technique \cite[Sec.~2]{action_3} for computing solutions for non-monotone systems of nonlinear algebraic equations that are not, in the first instance, the gradients of a scalar objective, as can arise in the local material update of classical plasticity models. In these situations, the dual scheme always produces symmetric Jacobians and, as is well-appreciated in computational plasticity circles, this is of practical significance. 

From the perspective of robust computation of approximate solutions of the scheme, the `universal' degenerate ellipticity of the scheme appears to make it particularly suitable for the application of Discontinuous Galerkin methods for elliptic problems \cite{arnold2000discontinuous}.

In closing, we mention that the ideas presented herein have strong links to modern mathematical thinking on Hidden Convexity in PDE advanced in \cite{brenier2018initial, brenier_book} (with the terminology of `Hidden Convexity' credited by Brenier to L.~C.~Evans), and appear to be also related to the recent work \cite{rockafellar2023augmented} on Hidden Convexity in Augmented Lagrangian techniques.

\appendix 
\renewcommand{\thesection}{\Alph{section}}

\section*{Appendices}

\section{Simple examples of the formalism and a conceptual outline}\label{app:rev_duality}

\subsection{An optimization problem for an algebraic system of equations}\label{sec:fin_dim}
 This section is excerpted from \cite{action_3}.
 
Consider a generally nonlinear system of algebraic equations in the variables $x \in \mathbb{R}^n$ given by
	\begin{equation}\label{eq:alg_sys}
	    A_\alpha (x) = 0,
	\end{equation}
	where $A: \mathbb{R}^n \to \mathbb{R}^N$ is a given function (a simple example would be $A_\alpha (x) = \bar{A}_{\alpha i} \, x^i - b_\alpha$, $\alpha = 1 \ldots N, i = 1 \ldots n$, where $\bar{A}$ is a constant matrix, \textit{not necessarily symmetric} (when $n = N$), and $b$ is a constant vector). We allow for all possibilities $0 < n \lesseqqgtr N > 0$.
	
	The goal is to construct an objective function whose critical points solve the system \eqref{eq:alg_sys} (when a solution exists) by defining an appropriate $x^* \in \mathbb{R}^n$ satisfying  $A_\alpha (x^*) = 0$.
	
	For this, consider first the auxiliary function
	\begin{equation*}
	    \widehat{S}_H(x,z) = z^\alpha A_\alpha (x) + H(x)
	\end{equation*}
	(where $H$ belongs to a class of scalar-valued function to be defined shortly) and define
	\begin{equation*}
	    S_H(z) = z^\alpha A_\alpha(x_H (z)) + H(x_H(z))
	\end{equation*}
	with the requirement that the system of equations
	\begin{equation}\label{eq:H_fin_dim}
	    z^\alpha \frac{\p A_\alpha}{\p x^i}(x) + \frac{\p H}{\p x^i}(x) = 0
	\end{equation}
	be solvable for the function $x = x_H(z)$ through the choice of $H$, and \textit{any} function $H$ that facilitates such a solution qualifies for the proposed scheme. 
	
	In other words, given a specific $H$, it should be possible to define a function $x_H(z)$ that satisfies 
	\begin{equation*} 
	z^\alpha \p_{x^i} A_\alpha (x_H(z)) + \p_{x^i} H (x_H(z)) = 0 \quad \forall z \in \mathbb{R}^N
	\end{equation*}
	(the domain of the function $x_H$ may accommodate more intricacies, but for now we stick to the simplest possibility). Note that \eqref{eq:H_fin_dim} is a set of $n$ equations in $n$ unknowns regardless of $N$ ($z$ for this argument is a parameter).
	
	Assuming this is possible, we have
	\begin{equation*}
	    \frac{\p S_H}{\p z^\beta} (z) = A_\beta(x_H(z)) +  \left( z^\alpha \frac{\p A_\alpha}{\p x^i}(x_H(z)) + \frac{\p H}{\p x^i}(x_H(z)) \right) \frac{\p x^i_H}{\p z^\beta}(z) = A_\beta(x_H(z)),
	\end{equation*}
	using \eqref{eq:H_fin_dim}. Thus,
	\begin{itemize}
	    \item if $z_0$ is a critical point of the objective function $S_H$ satisfying $\p_{z^\beta} S_H(z_0) = 0$, then the system $A_\alpha(x) = 0$ has a solution defined by $x_H(z_0)$; 
	    \item if the system $A_\alpha(x) = 0$ has a unique solution, say $y$, and if $z^H_0$ is any critical point of $S_H$, then $x_H\left(z^H_0 \right) = y$, for all admissible $H$.
	    \item If $A_\alpha(x) = 0$ has non-unique solutions, but $\p_{z^\beta} S(z) = 0$ ($N$ equations in $N$ unknowns) has a unique solution for a specific choice of the function $z \mapsto x_H(z)$ related to a choice of $H$, then such a choice of $H$ may be considered a selection criterion for imparting uniqueness to the problem $A_\alpha(x) = 0$.
        \item Finally, to see the difference of this approach with the Least-Squares (LS) Method, we note that the optimality condition for the objective $A_\alpha(x) A_\alpha(x)$ is $A_\alpha(x) \p_{x^i} A_\alpha(x) = 0 \centernot \implies A_\alpha(x) = 0$. 
        
        For a linear system $\bar{A} x = b$, the LS governing equations are given by
        \[
        \bar{A}^T \bar{A} z = \bar{A}^T b,
        \]
        with LS solution defined as $z$ even when the original problem $\bar{A} x = b$ does not have a solution (i.e., when $b$ is not in the column space of $\bar{A}$). The LS problem always has a solution, of course. In contrast, in the present duality-based approach with quadratic $H(x) = \frac{1}{2} x^T x$ the governing equation is
        \[
        - \bar{A}\bar{A}^T z = b
        \]
        with solution to $\bar{A} x = b$ given by $x = - \bar{A}^T z$, and the problem has a solution only when $\bar{A} x = b$ has a solution, since the column spaces of the matrices $\bar{A}$ and $\bar{A}\bar{A}^T$ are identical.
        \end{itemize}
\subsection{Application of base states: a simple example}\label{app:base_state_example}

 This section is excerpted from \cite{KA2}.
\begin{figure}
    \centering
    \includegraphics[width=0.3\textwidth]{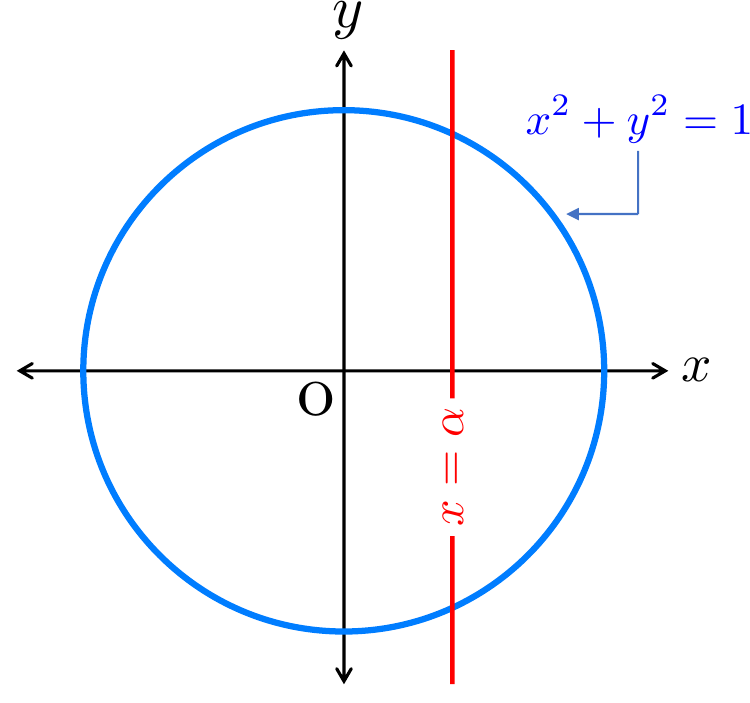}
    \caption{Schematic of Circle-line intersection.}
    \label{fig:circle_line}
\end{figure}

Consider the following algebraic system of equations for $(x,y)\in\mathbb{R}^2, \alpha \in \mathbb{R}$:
\begin{equation}
\begin{aligned}
    x^2 + y^2 & = 1\\
    x & = \alpha.
\end{aligned}
\label{eq:circle_line_primal}%
\end{equation}
A schematic is shown in Fig.~\ref{fig:circle_line}, and  solutions are given by $\{(\alpha, \pm \sqrt{1-\alpha^2}): |\alpha| \leq 1\}$.
We use the logic of Appendix \ref{sec:fin_dim} to multiply each of the above equations by a dual multiplier and add a quadratic auxiliary potential $H$:
\begin{equation*}
    \widehat{S}(x,y,\lambda,\gamma) = \lambda (x^2 + y^2 -1) + \gamma (x - \alpha) + \frac{1}{2}(x - \bar{x})^2 + \frac{1}{2}(y - \bar{y})^2,
\end{equation*}
where $\bar{x}, \bar{y}\in \mathbb{R}$ are constants (in this algebraic problem), and we will refer to the pair as a `base state.' Next we need to generate the analog of the mapping function $x_H$ of Appendix \ref{sec:fin_dim}, also referred to as the DtP mapping in the text:
\begin{equation} \label{eq:line_circle_dtp}
    \begin{aligned}
        \frac{\p S}{\p x} = 0 : \quad 2 \lambda x + \gamma + (x-\bar{x}) = 0 \quad & \Rightarrow \quad \mbox{for} \ \lambda \neq - \half, \quad x_H(\lambda,\gamma) = \frac{\bar{x} - \gamma}{2 \lambda + 1}; \\
        \frac{\p S}{\p y} = 0 : \quad 2 \lambda y + (y-\bar{y}) = 0 \quad & \Rightarrow \quad \mbox{for} \ \lambda \neq - \half, \quad  y_H(\lambda) = \frac{\bar{y}}{2 \lambda + 1}.
    \end{aligned} 
\end{equation}
{\boldmath$\lambda \neq - \half$}: Considering only the case $\lambda \neq - \half$ for the moment, the dual objective function is now obtained by substituting the DtP mapping into $\widehat{S}$:
\begin{equation*}
    S(\lambda,\gamma) = \lambda \left(x_H^2(\lambda,\gamma) + y_H^2(\lambda,\gamma) -1 \right) + \gamma \left(x_H(\lambda,\gamma) - \alpha \right) + \frac{1}{2}\left(x_H(\lambda,\gamma) - \bar{x} \right)^2 + \frac{1}{2}\left(y_H(\lambda,\gamma) - \bar{y} \right)^2.
\end{equation*}
The critical point equations for this objective, by design, are the equations \eqref{eq:circle_line_primal} with the substitution $(x \to x_H, y \to y_H)$:
\begin{equation*}
\begin{aligned}
    \frac{(\bar{x}-\gamma)^2}{(2\lambda+1)^2}+\frac{\bar{y}^2}{(2\lambda+1)^2} & = 1 \\ 
    \frac{\bar{x} - \gamma}{2\lambda+1} & = \alpha.
\end{aligned}
\label{eq:circle_line_dual}%
\end{equation*}
A necessary condition for solutions is
\[
(2\lambda + 1)^2 (1 - \alpha^2) - \bar{y}^2  = 0
\] 
which implies that dual solutions exist only for $|\alpha| \leq 1$ and when $|\alpha| = 1$ only if $\bar{y} = 0$. 

Thus, for $|\alpha| < 1$, $\bar{y} \neq 0$,
\begin{equation*}
    \lambda = \frac{1}{2}\left(\pm\frac{|\bar{y}|}{\sqrt{1-\alpha^2}}-1\right); \qquad \gamma = \bar{x} - \alpha(2\lambda + 1).
\end{equation*}
are the extrema of $S$.

Since only $\lambda \neq - \half$ is being considered, we do not consider the case $|\alpha| < 1$, $\bar{y} = 0$ here.

For $|\alpha| = 1$, $\bar{y} = 0$, the pairs
\begin{equation*}
    - \half \neq \lambda \in \R \ \mbox{arbitrary}; \qquad \gamma = \bar{x} - \alpha(2\lambda + 1).
\end{equation*}
are the extrema of $S$. 

Putting these dual solutions back into the DtP mapping \eqref{eq:line_circle_dtp}, we recover the correct primal solutions as expected.

\noindent {\boldmath$\lambda = - \half$}: When $\lambda = - \half$, extrema of $S$ exist only for $\bar{y} = 0$ and are given by $(\lambda = -\half, \gamma = \bar{x})$. However, note that this class of dual solutions does not define solutions to the primal problem in a non-vacuous manner (i.e.~while all primal solutions are admitted, it does not provide specific guidance for generating primal solutions).

The following conclusion can be drawn from this simple, yet non-trivial, example:
\begin{itemize}
    \item A `good' choice of the auxiliary function can be crucial for the success of the dual scheme in generating solutions to the primal problem. For example, if `no base states' are invoked within this class of quadratic auxiliary functions, i.e.~ $(\bar{x}, \bar{y}) = (0,0)$, while dual extrema exist, the scheme essentially fails to define primal solutions, except for the case $|\alpha| = 1$. And when the latter is the case for the primal problem, then a base state with non-zero $\bar{y}$ is not feasible for defining dual and primal solutions to the problem.
\end{itemize}
 \subsection{The idea behind the general formalism}
 This section is excerpted from \cite{action_2}.

The proposed scheme for generating variational principles for nonlinear PDE systems may be abstracted as follows: We first pose the given system of PDE as a \textit{first-order} system (introducing extra fields representing (higher-order) space and time derivatives of the fields of the given system); as before let us denote this collection of primal fields by $U$. `Multiplying' the primal equations by dual Lagrange multiplier fields, the collection denoted by $D$, adding a function $H(U)$, solely in the variables $U$ (the purpose of which, and associated requirements, will be clear shortly), and integrating by parts over the space-time domain, we form a `mixed' functional in the primal and dual fields given by
\begin{equation*}
    \widehat{S}_H [U,D] = \int_{[0,T]\times \Omega} \ \scl_H (\dee,U) \, dt dx 
\end{equation*}
where $\dee$ is a collection of local objects in $D$ and at most its first order derivatives. We then require that the family of functions $H$ be such that it allows the definition of a function $U_H(\dee)$ such that
\begin{equation*}
    \frac{\p \scl_H}{\p U} (\dee, U_H(\dee)) = 0
\end{equation*}
so that the \textit{dual} functional, defined solely on the space of the dual fields $D$, given by
\begin{equation*}
    S_H[D] = \int_{[0,T]\times \Omega} \scl_H(\dee,U_H(\dee)) \, dt dx 
\end{equation*}
has the first variation
\begin{equation*}
    \delta S_H = \int_{[0,T]\times \Omega} \frac{\p \scl_H}{\p \dee} \delta \dee \, dt dx.
\end{equation*}
By the process of formation of the functional $\widehat{S}_H$, it can then be seen that the (formal) E-L equations arising from $\delta S_H$ have to be the original first-order primal system, with $U$ substituted by $U_H(\dee)$, regardless of the $H$ employed.

Thus, the proposed scheme may be summarized as follows: we wish to pursue the following (local-global) critical point problem
\begin{equation*}
   \begin{smallmatrix} \mbox{extremize}\\ D\end{smallmatrix} \int_{[0,T]\times \Omega} \begin{smallmatrix} \mbox{extremize}\\ U\end{smallmatrix} \  \scl_H (\dee(t,x),U) \, dt dx,
\end{equation*}
where the pointwise extremization of $\scl_H$ over $U$, for fixed $\dee$, is made possible by the choice of $H$.

Furthermore, assume the Lagrangian $\scl_H$ can be expressed in the form
\begin{equation*}
    \scl_H(\dee, U) := - P(\dee)\cdot U + f(U,D) + H(U)
\end{equation*}
for some function $P$ defined by the structure of the primal first-order system ((linear terms in) first derivatives of $U$ after multiplication by the dual fields and integration by parts always produce such terms), and for some function $f$ which, when non-zero, does not contain any linear dependence in $U$. Our scheme requires the existence of a function $U_H$ defined from `solving $\frac{\p \scl}{\p U} (\dee, U) = 0$ for $U$,' i.e.~$\exists \  U_H(P(\dee),\dee)$ s.t. the equation
\begin{equation*}
    - P(\dee) + \frac{\p f}{\p U}(U_H(P(\dee),\dee), \dee) + \frac{\p H}{\p U}\left(U_H(P(\dee),\dee)\right) = 0
\end{equation*}
is satisfied. This requirement may be understood as follows: define
\begin{equation*}
    f(U, \dee) + H(U) =: M(U, \dee)
\end{equation*}
and assume that it is possible, through the choice of $H$, to make the function $\frac{\p M}{\p U}(U, \dee)$ \textit{monotone} in $U$ so that a function $U_H(P,\dee)$ can be defined that satisfies
\begin{equation*}
    \frac{\p M}{\p U}(U_H(P,\dee), \dee) = P, \quad \forall P.
\end{equation*}
Then the Lagrangian is
\begin{equation*}
    \scl(\dee, U_H(P(\dee),\dee)) = - P(\dee) \cdot U_H(P(\dee),\dee) + M(U_H(P(\dee),\dee), \dee) =: - M^*(P(\dee), \dee)
\end{equation*}
where $M^*(P,\dee)$ is the Legendre transform of the function $M$ w.r.t $U$, with $\dee$ considered as a parameter.

Thus, our scheme may also be interpreted as designing a concrete realization of abstract saddle point problems in optimization theory \cite{rockafellar1974conjugate}, where we exploit the fact that, in the context of `solving' PDE viewed as constraints implemented by Lagrange multipliers to generate an unconstrained problem, there is a good deal of freedom in choosing an objective function(al) to be minimized. We exploit this freedom in choosing the function $H$ to develop dual variational principles corresponding to general systems of PDE.

\section{Dirichlet boundary conditions for elastostatics in first-order form \eqref{eq:gov_cont_mech}}\label{app:A}
Consider the system \eqref{eq:elastostat} with 
\[
U = (y_1, y_2, y_3, F_{11}, F_{12},F_{13}, F_{21}, F_{22}, F_{23}, F_{31}, F_{32}, F_{33}).
\]
For $\Gamma = 4, \dots, 12; j = 1, \ldots, 3$, we consider $F_{\Gamma j}$ to be of the form
\[
F_{\Gamma j} (U) := \calA_{\Gamma Ij} U_I
\]
with $\calA, \calB$  constant matrices with $0$ entries, unless otherwise specified. Then,
\begin{equation*}
    \begin{aligned}
        \calA_{411} = 1; & \qquad \calB_{44} = 1 & \qquad \Longrightarrow \qquad \p_1 y_1 - F_{11} = 0\\
        \calA_{512} = 1; & \qquad \calB_{55} = 1 & \qquad \Longrightarrow \qquad \p_2 y_1 - F_{12} = 0\\
        \calA_{613} = 1; & \qquad \calB_{66} = 1 & \qquad \Longrightarrow \qquad \p_3 y_1 - F_{13} = 0 \\
        \calA_{721} = 1; & \qquad \calB_{77} = 1 & \qquad \Longrightarrow \qquad \p_1 y_2 - F_{21} = 0\\
        \calA_{822} = 1; & \qquad \calB_{88} = 1 & \qquad \Longrightarrow \qquad \p_2 y_2 - F_{22} = 0\\
        \calA_{923} = 1; & \qquad \calB_{99} = 1 & \qquad \Longrightarrow \qquad \p_3 y_2 - F_{23} = 0 \\
        \calA_{10 \, 31} = 1; & \qquad \calB_{10 \, 10} = 1 & \qquad \Longrightarrow \qquad \p_1 y_3 - F_{31} = 0\\
        \calA_{11 \, 32} = 1; & \qquad \calB_{11 \, 11} = 1 & \qquad \Longrightarrow \qquad \p_2 y_3 - F_{32} = 0\\
        \calA_{12 \, 33} = 1; & \qquad \calB_{12 \, 12} = 1 & \qquad \Longrightarrow \qquad \p_3 y_3 - F_{33} = 0.\\
    \end{aligned}
\end{equation*}

Let the matrix entries $B_{\Gamma j} = 0$ unless otherwise specified and $n$ be the outward unit normal field on the boundary $\p \Omega$.

Now, let $y^*_1$ be the desired Dirichlet b.c.~on $y_1$ on $\p \Omega_4 = \p \Omega_5 = \p \Omega_6 =: \p \Omega_{456}$, and for $\Gamma = 4,5,6$ let $B_{\Gamma j} n_j$ be defined as 
\[
B_{\Gamma j} n_j := \calA_{\Gamma I j} n_j y^*_I \qquad \mbox{on} \ \p \Omega_{456},
\]
with $y^*_I = 0, I \neq 1$ without loss of generality.
Then\eqref{eq:gov_bc} implies the Dirichlet b.c.
\[
(y_1 - y_1^*) n_j = 0 \ \forall \ j = 1,2,3 \ \mbox{on} \  \Longrightarrow  y_1 - y_1^* = 0  \ \mbox{on} \ \p \Omega_{456}.
\]
Similarly, let $y^*_2$ be the desired Dirichlet b.c.~on $y_2$ on $\p \Omega_7 = \p \Omega_8 = \p \Omega_9 =: \p \Omega_{789}$, and for $\Gamma = 7,8,9$ let $B_{\Gamma j} n_j$ be defined as 
\[
B_{\Gamma j} n_j := \calA_{\Gamma I j} n_j y^*_I \qquad \mbox{on} \ \p \Omega_{789},
\]
with $y^*_I = 0, I \neq 2$ without loss of generality.
Then\eqref{eq:gov_bc} implies the Dirichlet b.c.
\[
(y_2 - y_2^*) n_j = 0 \ \forall \ j = 1,2,3 \ \mbox{on} \  \Longrightarrow  y_2 - y_2^* = 0  \ \mbox{on} \ \p \Omega_{789},
\]
and let $y^*_3$ be the desired Dirichlet b.c.~on $y_3$ on $\p \Omega_{10} = \p \Omega_{11} = \p \Omega_{12} =: \p \Omega_{10 \, 11 \, 12}$, and for $\Gamma = 10, 11, 12$ let $B_{\Gamma j} n_j$ be defined as 
\[
B_{\Gamma j} n_j := \calA_{\Gamma I j} n_j y^*_I \qquad \mbox{on} \ \p \Omega_{10 \, 11 \, 12},
\]
with $y^*_I = 0, I \neq 3$ without loss of generality.
Then\eqref{eq:gov_bc} implies the Dirichlet b.c.
\[
(y_3 - y_3^*) n_j = 0 \ \forall \ j = 1,2,3 \ \mbox{on} \  \Longrightarrow  y_3 - y_3^* = 0  \ \mbox{on} \ \p \Omega_{10 \, 11 \, 12}.
\]

\bibliographystyle{alphaurl}\bibliography{dual_plasticity.bib}
\end{document}